\documentclass[12pt,a4paper]{amsart}
\usepackage{graphics}
\usepackage{epsfig}
\usepackage{graphicx}
\theoremstyle{plain}
\usepackage{amssymb}

\advance\hoffset-20mm \advance\textwidth40mm

\newtheorem{theorem}{Theorem}

\newtheorem*{theo*}{Theorem}

\theoremstyle{definition}

\newtheorem*{definition*}{Definition}

\def\AA{{\mathbb A}}
\def\PP{{\mathbb P}}
\def\ZZ{{\mathbb Z}}
\def\GG{{\mathbb G}}
\def\Of{{\mathcal{O}}}
\def\fg{{\mathfrak g}}
\def\fh{{\mathfrak h}}
\def\ft{{\mathfrak t}}
\def\fp{{\mathfrak q}}
\def\Aut{\mathop{\rm Aut}}
\def\Lie{\mathop{\rm Lie}}
\def\deg{\mathop{\rm deg}}
\def\SL{\mathop{\rm SL}}
\def\Sp{\mathop{\rm Sp}}
\def\SO{\mathop{\rm SO}}
\def\PSL{\mathop{\rm PSL}}
\def\PSp{\mathop{\rm PSp}}
\def\PSO{\mathop{\rm PSO}}

\begin{document}
\sloppy
\title[Flag varieties as compactifications of $\GG_a^n$]
{Flag varieties as equivariant compactifications of $\GG_a^n$}
\author[I.V. Arzhantsev]{Ivan V. Arzhantsev}
\thanks{Supported by RFBR grants 09-01-00648-a, 09-01-90416-Ukr-f-a, and the Deligne
2004 Balzan prize in mathematics.}
\address{Department of Algebra, Faculty of Mechanics and Mathematics,
Moscow State University, Leninskie Gory 1, GSP-1, Moscow, 119991,
Russia } \email{arjantse@mccme.ru}
\date{\today}
\begin{abstract}
Let $G$ be a semisimple affine algebraic group and $P$ a parabolic subgroup of $G$.
We classify all flag varieties $G/P$ which admit an action of the commutative
unipotent group $\GG_a^n$ with an open orbit.
\end{abstract}
\maketitle
\section*{Introduction}

Let $G$ be a connected semisimple affine algebraic group of adjoint type over
an algebraically closed field
of characteristic zero, and $P$ be a parabolic subgroup of $G$. The
homogeneous space $G/P$ is called a (generalized) flag variety. Recall
that $G/P$ is complete and the action of the unipotent radical $P_u^{-}$ of the opposite parabolic subgroup $P^-$
on $G/P$ by left multiplication is generically transitive. The open orbit $\Of$ of this action is called the big Schubert cell on $G/P$. Since $\Of$ is isomorphic to the affine space $\AA^n$, where $n=\dim G/P$, every flag variety may be regarded as a compactification of an affine space.

Notice that the affine space $\AA^n$ has a structure of the vector group, or,
equivalently, of the commutative unipotent affine algebraic group $\GG_a^n$. We say
that a complete variety $X$ of dimension $n$ is an equivariant compactification of the group $\GG_a^n$, if there exists a regular action $\GG_a^n\times X \to X$ with a dense open orbit.
A systematic study of equivariant compactifications of the group $\GG_a^n$ was initiated
by B.~Hassett and Yu.~Tschinkel in~\cite{HT}, see also \cite{Sh} and \cite{ASh}.

In this note we address the question whether a flag variety $G/P$ may be realized as an
equivariant compactification of $\GG_a^n$. Clearly, this is the case when the group $P_u^{-}$,
or, equivalently, the group $P_u$ is commutative.
It is a classical result that the connected component
$\widetilde{G}$ of the automorphism group of the variety $G/P$ is a semisimple group of adjoint type, and $G/P=\widetilde{G}/Q$ for some parabolic subgroup $Q\subset \widetilde{G}$.
In most cases the group $\widetilde{G}$ coincides with $G$,
and all exceptions are well known, see \cite{On1}, \cite[Theorem~7.1]{On2},
\cite[page~118]{Ti}, \cite[Section~2]{De}.
If $\widetilde{G} \ne G$, we say that
$(\widetilde{G}, Q)$ is the covering
pair of the exceptional pair $(G,P)$. For a simple group $G$, the exceptional pairs are
$(\PSp(2r), P_1)$, $(\SO(2r+1), P_r)$ and $(G_2, P_1)$ with the covering pairs
$(\PSL(2r), P_1)$, $(\PSO(2r+2), P_{r+1})$ and $(\SO(7), P_1)$ respectively,
where $PH$ denotes the quotient of the group $H$ by its center, and $P_i$ is the maximal parabolic subgroup associated with the $i$th simple root. It turns out that for a simple group $G$
the condition $\widetilde{G}\ne G$ implies that
the unipotent radical $Q_u$ is commutative
and $P_u$ is not. In particular, in this case $G/P$ is an equivariant compactification of $\GG_a^n$. Our main result states that these are the only possible cases.

\begin{theorem} \label{main}
Let $G$ be a connected semisimple group of adjoint type and $P$ a parabolic subgroup of $G$. Then
the flag variety $G/P$ is an equivariant compactification of $\GG_a^n$ if and only if
for every pair $(G^{(i)}, P^{(i)})$, where $G^{(i)}$ is a simple component of $G$
and $P^{(i)}= G^{(i)} \cap P$, one of the following conditions holds:
\begin{enumerate}
\item
the unipotent radical $P^{(i)}_u$ is commutative;
\item
the pair $(G^{(i)},P^{(i)})$ is exceptional.
\end{enumerate}
\end{theorem}

For convenience of the reader, we list all pairs $(G,P)$, where $G$ is a simple
group (up to local isomorphism) and $P$ is a parabolic subgroup with a commutative unipotent radical:
$$
(\SL(r+1),\, P_i), \ i=1,\ldots,r; \quad (\SO(2r+1),\, P_1); \quad (\Sp(2r),\, P_r);
$$
$$
(\SO(2r),\, P_i), \ i=1,r-1,r; \quad (E_6,\, P_i), \ i=1,6; \quad (E_7,\, P_7),
$$
see~\cite[Section~2]{RRS}. The simple roots $\{\alpha_1,\ldots,\alpha_r\}$ are indexed
as in \cite[Planches I-IX]{Bu}. Note that the unipotent radical of $P_i$ is commutative
if and only if the simple root $\alpha_i$ occurs in the highest root $\rho$ with
coefficient 1, see \cite[Lemma~2.2]{RRS}.
Another equivalent condition is that the fundamental
weight $\omega_i$ of the dual group $G^\vee$ is minuscule, i.e., the weight system of the simple $G^\vee$-module $V(\omega_i)$ with the highest weight $\omega_i$ coincides with the orbit $W\omega_i$ of the Weyl group $W$.

\section{Proof of Theorem~\ref{main}}

If the unipotent radical $P_u^-$ is commutative, then the action of $P_u^-$ on $G/P$
by left multiplication is the desired generically transitive $\GG_a^n$-action, see,
for example, \cite[pp.~22-24]{LR}. The same arguments work when for the connected component
$\widetilde{G}$ of the automorphism group $\Aut(G/P)$ one has $G/P=\widetilde{G}/Q$
and the unipotent radical $Q_u^-$ is commutative. Since
$$
G/P \, \cong \, G^{(1)}/P^{(1)} \times \ldots \times G^{(k)}/P^{(k)},
$$
where $G^{(1)},\ldots,G^{(k)}$ are the simple components of the group $G$,
the group $\widetilde{G}$ is isomorphic to the direct product
$\widetilde{G^{(1)}}\times \ldots \times \widetilde{G^{(k)}}$, cf.~\cite[Chapter~4]{On3}.
Moreover,
$Q_u \cong Q_u^{(1)}\times \ldots \times Q_u^{(k)}$ with $Q^{(i)}=\widetilde{G^{(i)}}\cap Q$,
Thus the group $Q_u^-$
is commutative if and only if for every pair $(G^{(i)}, P^{(i)})$ either $P_u^{(i)}$
is commutative or the pair $(G^{(i)}, P^{(i)})$ is exceptional.

Conversely, assume that $G/P$ admits a generically transitive $\GG_a^n$-action.
One may identify $\GG_a^n$ with a commutative unipotent subgroup $H$ of $\widetilde{G}$,
and the flag variety $G/P$ with $\widetilde{G}/Q$, where $Q$ is a parabolic subgroup
of $\widetilde{G}$.

Let $T\subset B$ be a maximal torus and a Borel subgroup of the group $\widetilde{G}$
such that $B\subseteq Q$.
Consider the root system $\Phi$ of the tangent algebra $\fg=\Lie(\widetilde{G})$ defined by the
torus $T$, its decomposition $\Phi=\Phi^+ \cup \Phi^-$ into positive and negative roots associated with $B$, the set of simple roots $\Delta \subseteq \Phi^+$,
$\Delta=\{\alpha_1,\ldots,\alpha_r\}$, and the root decomposition
$$
\fg \ = \ \bigoplus_{\beta\in\Phi^-} \fg_{\beta} \ \oplus \ft \ \oplus
\bigoplus_{\beta\in\Phi^+} \fg_{\beta},
$$
where $\ft=\Lie(T)$ is a Cartan subalgebra in $\fg$ and $\fg_{\beta}=\{x\in\fg : [y,x]=\beta(y)x
\ \text{for all} \ y\in\ft\}$ is the root subspace.
Set $\fp=\Lie(Q)$ and $\Delta_Q=\{\alpha\in\Delta : \fg_{-\alpha} \nsubseteq \fp\}$.
For every root $\beta = a_1\alpha_1 + \ldots + a_r\alpha_r$ define
$\deg(\beta)=\sum_{\alpha_i\in \Delta_P} a_i$. This gives a $\ZZ$-grading on the Lie algebra~$\fg:$
$$
\fg = \bigoplus_{k\in\ZZ} \fg_k, \quad \text{where} \quad \ft\subseteq\fg_0 \quad
\text{and} \quad \fg_{\beta} \subseteq \fg_k \quad \text{with} \quad \ k=\deg(\beta).
$$
In particular,
$$
\fp \ = \ \bigoplus_{k\ge 0} \fg_k \quad \text{and} \quad
\fp_u^- \ = \ \bigoplus_{k<0} \fg_k.
$$
Assume that the unipotent radical $Q_u^-$ is not commutative, and consider
$\fg_{\beta} \subseteq [\fp_u^-,\fp_u^-]$. For every $x\in \fg_{\beta} \setminus \{0\}$
there exist $z'\in\fg_{\beta'}\subseteq\fp_u^-$ and $z''\in\fg_{\beta''}\subseteq\fp_u^-$
such that $x=[z',z'']$. In this case $\deg(z')>\deg(x)$ and $\deg(z'')>\deg(x)$.

Since the subgroup $H$ acts on $\widetilde{G}/Q$ with an open orbit, one may conjugate $H$
and assume that the $H$-orbit of the point $eQ$ is open in $\widetilde{G}/Q$. This implies
$\fg = \fp \oplus \fh$, where $\fh=\Lie(H)$. On the other hand, $\fg = \fp \oplus \fp_u^-$.
So every element $y\in\fh$ may be (uniquely) written as $y=y_1+y_2$, where $y_1 \in \fp$,
$y_2 \in \fp_u^-$, and the linear map $\fh \to \fp_u^-$, $y \mapsto y_2$, is bijective.
Take the elements $y,y',y'' \in \fh$ with $y_2=x,\, y'_2=z',\, y''_2=z''$.
Since the subgroup $H$ is commutative, one has $[y',y'']=0$. Thus
$$
[y'_1 + y'_2, y''_1+y''_2]\ =\ [y'_1,y''_1] + [y'_2,y''_1] + [y'_1,y''_2] + [y'_2,y''_2]\ =\ 0.
$$
But
$$
[y'_2,y''_2]=x \quad \text{and} \quad [y'_2,y''_1] + [y'_1,y''_2] + [y'_2,y''_2] \in
\bigoplus_{k>\deg(x)} \fg_k.
$$
This contradiction shows that the group $Q_u^-$ is commutative. As we have seen, the latter
condition means that for every pair $(G^{(i)}, P^{(i)})$ either
the unipotent radical $P^{(i)}_u$ is commutative or the pair $(G^{(i)},P^{(i)})$ is exceptional.
The proof of Theorem~\ref{main} is completed.

\section{Concluding remarks}

If a flag variety $G/P$ is an equivariant compactification of $\GG_a^n$, then
it is natural to ask for a classification of all generically transitive
$\GG_a^n$-actions on $G/P$ up to equivariant isomorphism. Consider the projective space
$\PP^n \cong \SL(n+1)/P_1$. In \cite{HT}, a correspondence between equivalence classes
of generically transitive $\GG_a^n$-actions on $\PP^n$ and isomorphism classes of
local (associative, commutative) algebras of dimension $n+1$ was established.
This correspondence together with classification results from \cite{ST} yields that
for $n\ge 6$ the number of equivalence classes of generically transitive $\GG_a^n$-actions
on $\PP^n$ is infinite, see \cite[Section~3]{HT}. On the contrary, a generically transitive
$\GG_a^n$-action on the non-degenerate projective quadric \linebreak $Q_n~\cong~\SO(n+2)/P_1$ is unique
\cite[Theorem~4]{Sh}. It would be interesting to study the same problem for the Grassmannians
$\text{Gr}(k,r+1) \cong \SL(r+1)/P_k$, where $2\le k\le r-1$.

\section*{Acknowledgement}

The author is indebted to N.A.~Vavilov for a discussion which results in this note.
Thanks are also due to D.A.~Timashev and M.~Zaidenberg for their interest and valuable comments.


%

\begin{thebibliography}{}
%
\bibitem{ASh}
I.V.~Arzhantsev and E.V.~Sharoyko,
Hassett-Tschinkel correspondence: modality and projective hypersurfaces.
arXiv:0912.1474 [math.AG]
%
\bibitem{Bu}
N.~Bourbaki, Groupes et alg\'ebres de Lie, Chaps. 4,5 et 6. Paris, Hermann, 1975.
%
\bibitem{De}
M.~Demazure,
Automorphismes et d\'eformations des vari\'et\'es de Borel.
Invent. Math. 39 (1977), 179--186.
%
\bibitem{HT}
B.~Hassett and Yu.~Tschinkel,
Geometry of equivariant compactifications of $\Bbb{G}_a^n$.
Int. Math. Res. Notices~22 (1999), 1211--1230.
%
\bibitem{LR}
V.~Lakshmibai and K.N.~Raghavan, Standard Monomial Theory.
Invariant Theoretic Approach. Encyclopaedia of Mathematical
Sciences, Vol.~137, Springer, 2008.
%
\bibitem{On1}
A.L.~Onishchik,
On compact Lie groups transitive on certain manifolds.
Dokl. Akad. Nauk SSSR 135 (1961), 531--534 (Russian);
English transl.: Sov. Math., Dokl. 1 (1961), 1288--1291.
%
\bibitem{On2}
A.L.~Onishchik,
Inclusion relations between transitive compact transformation groups.
Tr. Mosk. Mat. O.-va~11 (1962), 199--242 (Russian).
%
\bibitem{On3}
A.L.~Onishchik,
Topology of transitive transformation groups.
Leipzig: Johann Ambrosius Barth., 1994. 
%
\bibitem{RRS}
R.~Richardson, G.~R\"ohrle and R.~Steinberg,
Parabolic subgroups with Abelian unipoten radical.
Invent. Math 110 (1992), 649--671.
%
\bibitem{Sh}
E.V.~Sharoyko,
Hassett-Tschinkel correspondence and automorphisms of the quadric.
Sbornik Math. 200 (2009), no.~11, 145--160.
%
\bibitem{ST}
D.A.~Suprunenko and R.I.~Tyshkevich,
Commutative matrices. Academic Press, New York, 1969.
%
\bibitem{Ti}
J.~Tits,
Espaces homog\'enes complexes compacts.
Comm. Math. Helv. 37 (1962), 111--120.
%
\end{thebibliography}
\end{document}